\documentclass[twoside,reqno]{amsart}

\usepackage{latexsym}

\newcommand{\binm}{\binom}

\newcommand{\be}{\begin{equation}}
\newcommand{\ee}{\end{equation}}
\newcommand{\ba}{\begin{array}}
\newcommand{\ea}{\end{array}}
\newcommand{\bmn}{\begin{eqnarray}}
\newcommand{\emn}{\end{eqnarray}}
\newcommand{\bnm}{\begin{eqnarray*}}
\newcommand{\enm}{\end{eqnarray*}}
\newcommand{\bln}{\begin{subequations}}
\newcommand{\eln}{\end{subequations}}

\newtheorem{entry}{Entry}

\newcommand{\bbtm}[4]{\bibitem{kn:#1}{#2,}~{#3,}~{#4.}}
\newcommand{\cito}[1]{\cite{kn:#1}}
\newcommand{\citu}[2]{\cite[#2]{kn:#1}}

\begin{document} 
{} 
\title{Verification of Binomial theorem and Chu-Vandermonde convolution by \\the finite difference method}
\author{$^a$Chuanan Wei, $^b$Dianxuan Gong}
\dedicatory{
$^A$Department of Information Technology\\
  Hainan Medical College,  Haikou 571101, China\\
         $^B$College of Sciences\\
             Hebei Polytechnic University, Tangshan 063009, China}
\thanks{Corresponding author$^*$. \emph{Email address}:
      gongdianxuan@yahoo.com.cn}

\address{ }


\begin{abstract}
In this note, we show that Binomial theorem and Chu-Vandermonde
convolution can both be verified by the finite difference method.
\end{abstract}

\maketitle\thispagestyle{empty}
\markboth{Chuanan Wei, Dianxuan Gong}
         {Verification of Binomial theorem and Chu-Vandermonde convolution}


There are numerous binomial coefficient identities in the
literature. Thereinto, Binomial theorem and Chu-Vandermonde
convolution(cf. Bailey \citu{bailey}{$\xi$1.3}) can be stated,
respectively, as
 \bmn
 &&\sum_{k=0}^n\binm{n}{k}x^k=(1+x)^n, \label{binomial}\\
&&\sum_{k=0}^n\binm{x}{k}\binm{y}{n-k}=\binm{x+y}{n}.
\label{vandermonde}
 \emn

Let $\bigtriangleup$ be the usual difference operator with the unit
increment. For a complex function $f(\tau)$, the finite difference
of order $n$ can be calculated through the following Newton-Gregory
formula(cf. Graham et al. \citu{graham}{$\xi$5.3}):
 \bmn
\bigtriangleup^nf(\tau)=\sum_{k=0}^n(-1)^{n+k}\binm{n}{k}f(\tau+k).
\label{newton}
 \emn
When $f(\tau)$ is a polynomial of degree $m\leq n$, then
$\bigtriangleup^n f(\tau)$ vanishes for $0\leq m<n$ and otherwise,
equals $m!$ times the leading coefficient of $f(\tau)$ for $m=n$.
Based on \eqref{binomial} and  \eqref{vandermonde}, Chu~\cito{chu}
 derived respectively Abel's identities and Hagen-Rothe convolutions  by
the finite difference method. Inspired by the work just mentioned,
we shall show that \eqref{binomial} and \eqref{vandermonde} can also
be verified in the same way.

\textbf{Verification of Binomial theorem:}

  Define the function $f(\tau)$ by
\[f(\tau)=(-x)^{\tau}.\]
 On one hand, we obtain, according to
\eqref{newton}, the relation
 \bmn \label{bino-diff-a}
\bigtriangleup^nf(\tau)=\sum_{k=0}^n(-1)^{n+k}\binm{n}{k}(-x)^{\tau+k}.
 \emn
 On the other hand, it is not difficult to get the result
\bmn \label{bino-diff-b}
 \bigtriangleup^nf(\tau)=(-x)^{\tau}(-x-1)^n.
 \emn
Combing \eqref{bino-diff-a} with \eqref{bino-diff-b}, we derive the
equation
\[\sum_{k=0}^n(-1)^{n+k}\binm{n}{k}(-x)^{\tau+k}=(-x)^{\tau}(-x-1)^n.\]
Dividing both sides of the last equation by $(-1)^n(-x)^{\tau}$, we
deduce \eqref{binomial} to complete the verification.

\textbf{Verification of Chu-Vandermonde convolution:}

 It is easy to
 see that \eqref{vandermonde} is right when $n=0$. We shall also assume $n>0$ in the following
  argumentations. For a natural number
 $i$ with $1\leq i\leq n$, define the function $f(y)$ by
\[f(y)=\binm{y-i}{\,n-i\,}\]
which is a polynomial of degree $n-i<n$ in $y$. According to
\eqref{newton}, we get the relation
 \bmn \label{bino-diff-s}
\bigtriangleup^nf(y)=\sum_{k=0}^n(-1)^{n+k}\binm{n}{k}\binm{y+k-i}{n-i}=0.
 \emn
Define the function $F(x)$ by
\[F(x)=\sum_{k=0}^n\binm{x}{k}\binm{y}{n-k}\]
which is a polynomial of degree $n$ in $x$. Considering that
\[F(i-y-1)=\sum_{k=0}^n\binm{i-y-1}{k}\binm{y}{n-k}=\frac{\binm{y}{i}}{\binm{n}{i}}
\sum_{k=0}^n(-1)^{k}\binm{n}{k}\binm{y+k-i}{n-i},\] we asserts that
\eqref{bino-diff-s} leads to $F(i-y-1)=0$. Thereby, we find that all
the zeros of $F(x)$ are given by $\{i-y-1: 1\leq i\leq n\}$.
Observing further that the polynomial $\binm{x+y}{n}$ has the same
zeros as $F(x)$, there must exist a constant $\theta$ such that
\[F(x)=\theta\binm{x+y}{n}\quad\text{with}\quad \theta=1\]
where $\theta$ has been determined by letting $x=0$ in the last
equation. In conclusion, we have verified the correctness of
\eqref{vandermonde}.


\end{document}